\documentstyle[12pt]{article}

\def\fr{\frac}

\def\de{\delta}

\def\a{a^{\dagger}}

\newcommand{\bd}{\begin{displaymath}}
\newcommand{\ed}{\end{displaymath}}
\newcommand{\bb}{\begin{equation}}
\newcommand{\ee}{\end{equation}}

\begin{document}
\baselineskip 1.3 \baselineskip

%% Title Page 

\vspace{.2cm}

\begin{center}
\Large {\bf The Polar Decomposition of q-Deformed Boson Algebra}
\\[1cm]
\large W.-S.Chung \\[.3cm]
\normalsize  
Department of Physics and Research Institute of Natural Science, \\ 
\normalsize College of Natural Sciences,  \\
\normalsize  Gyeongsang National University,   \\
\normalsize   Jinju, 660-701, Korea
\end{center}

\vspace{0.5cm}
\begin{abstract}
In this paper we used the finite Fourier transformation to obtain the
polar decomposition of the q-deformed boson algebra with $q$ a root of
unity.

\end{abstract}

Since the q-deformed boson algebra was known [1,2], much work 
has been accomplished in this direction [3-8].
It is well known that this
algebra has an infinite dimensional representation when the
deformation parameter $q$ is real, while it has a finite dimensional
one when $q$ is a root of unity.

In this paper we use the finite Fourier transformation to obtain the
polar decomposition of the q-deformed boson algebra with $q$ a root
of unity. The polar decomposition of an operator means that the
operator can be written as a unitary times a hermitian operator.

Now let us consider the q-deformed boson algebra [1,2]. It is defined
as 
\bd
a\a-q\a a=q^{-N},
\ed
\bd
[N,\a]=\a,
\ed
\bb
[N,a]=-a.
\ee
Then, the step operators and the number operator are defined as
\bd
a=\sum_{n=0}^{\infty}\sqrt{[n]}|n-1><n|,
\ed
\bd
\a=\sum_{n=0}^{\infty}\sqrt{[n+1]}|n+1><n|,
\ed
\bb
N=\sum_{n=0}^{\infty}n|n><n|,
\ee
where the q-number $[x]$ is given by
\bd
[x]=\fr{q^x-q^{-x}}{q-q^{-1}}
\ed
If we take $q$ to be $q=e^{i\fr{2\pi}{s+1}}$, we have $[s+1]=0$, which
implies
\bb
\a|s>=0.
\ee
Thus, the spectrum becomes finite dimensional and the matrix
representation of the step and the number operators are given by
\bd
a=\sum_{n=0}^{s}\sqrt{[n]}|n-1><n|,
\ed
\bd
\a=\sum_{n=0}^{s}\sqrt{[n+1]}|n+1><n|,
\ed
\bb
N=\sum_{n=0}^{s}n|n><n|,
\ee
It can be easily checked that the step operators satisfy the
nilpotency condition;
\bb
a^{s+1}=(\a)^{s+1}=0.
\ee

\def\h{h^{\dagger}}

In this case the step operators can be decomposed into the
following two pices ;

\bd
a=\sqrt{\{g+1\}}\h=\h\sqrt{\{g\}},
\ed
\bb
\a=\sqrt{\{g\}}h=h\sqrt{\{g+1\}} ,
\ee
with $g=q^N$ and
\bd
h=\sum_{n=0}^{s-1}|n+1><n|=|1><0|+\cdots+|s><s-1|,
\ed
\bb
\h=\sum_{n=0}^{s-1}|n><n+1|=|0><1|+\cdots+|s-1><s|,
\ee
where we have used the symbols
\bd
\{g\}=\fr{g-g^{-1}}{q-q^{-1}}=\a a=[N],
\ed
\bb
\{g+1\}=\fr{qg-q^{-1}g^{-1}}{q-q^{-1}}=a\a =[N+1].
\ee     

Here we cannot call eq.(6) the polar decomposition because in the
finite dimensional representation with $q$ a root of unity
the operator $h$ and $\h$ do not satisfy the unitary condition.
 From the
definition of $h$ and $\h$ , the commutation relations of $(h,\h)$ and
$g$ become 
\bd 
gh=qhg,
\ed
\bb
g\h =q^{-1}\h g.
\ee
Since $a$ and $\a$ operator act on the basis of $(s+1)$-dimensional
Fock space $\{ |0>,\cdots,|s>\}$,  the opeartor $h$ and $\h$ are not
unitary any more. Instead $h$ and $\h$ satisfy 
\bd
h\h=1-|0><0|,
\ed
\bb
\h h =1-|s><s| .
\ee
Then it can be easily checked that the decomposition (6) satisfy the
relation 
\bb
\a a =[N],~~~~~a\a=[N+1].
\ee
From the definitin of $g$ and $ h$ with $q=e^{i\fr{2\pi}{s+1}}$, we have
\bb
g^{s+1}=1,~~~h^{s+1}=0.
\ee
Let us introduce the finite Fourier transformation 
\bb
F=\fr{1}{\sqrt{s+1}}\sum_{n,m=0}^sq^{mn}|m><m|
\ee
satisfying the unitary condition
\def\f{F^{\dagger}}
\bd
F\f=\f F=1.
\ed
Using this operator we have
\bd
h=Fg^{-1}\f -|0><s|,
\ed
\bb
\h=Fg\f -|s><0|.
\ee
It can be easily checked that the eq.(14) satisfies the relation (9).

\def\m{|\phi_m>}

\def\ma{|\phi_{m+1}>}

\def\mb{|\phi_{m-1}>}

The finite Fourier transformation dfines the phase state $\m=F|m>$ and
the phase state $\{|\phi_0>,\cdots,|\phi_s>\}$ constitute an
orthonormal set of states dual to that of the number states
$\{|0>,\cdots,|s>\}$: \bb
<\phi_m \m =\de_{mn}.
\ee
Further we can define creation and  annihilation operator acting in the
phase state
basis as follows;
\def\at{\tilde{a}}
\def\att{\tilde{a}^{\dagger}}
\bd
\at =Fa \f =\sum_{m=0}^s \sqrt{[m]}\mb <\phi_m|,
\ed                             
\bd
\att =F\a \f =\sum_{m=0}^s \sqrt{[m+1]}\ma <\phi_m|,
\ed                             
\def\n{\tilde{N}}
\bb
\n =F N \f =\sum_{m=0}^s m|\phi_m> <\phi_m|.
\ee
Then, we see that they satisfy
\bb
\at \att -q \att \at =q^{-\n}=H,
\ee
where 
\bd
H=h+|0><s|.
\ed
At this stage we notice that $H$ is a unitary operator although $h$ is
not  unitary. Thus, the new operators $(H,H^{\dagger})$ and $g$ obey
the following commutation relation \def\hh{H^{\dagger}}
\bd
gH=qHg,~~~~g\hh=q^{-1}\hh g,
\ed
\bb
H^{s+1}=1.
\ee
Therefore, we obtain the polar decomposition of $\at$ and $\att$,
which  is given by
\bd
\at=\sqrt{\{ \hh+1\}}g^{-1}=g^{-1}\sqrt{\{ \hh\}},
\ed
\bb
\att=\sqrt{\{ \hh\}}=g\sqrt{\{ \hh +1\}},
\ee
where
\bd
\{\hh\}=\fr{\hh-H}{q-q^{-1}},
\ed
\bb
\{\hh +1\}=\fr{q\hh-q^{-1}H}{q-q^{-1}}.
\ee
         
To conclude, we considered the finite dimensional representation of
the q-defomrmed boson algebra when $q$ is a root of unity. In this
case the step operators of this algebra do not seem to have the polar
decomposition ( a unitary times a hermitian operator ).
But we used the finite Fourier transformation and the unitary
transformation of the step and the number operators to obtain the
polar decomposition of the step opeartor of the q-deformed boson
algebra.
          
\section*{Acknowledgement}
This                   paper                was
supported         by  
the   KOSEF (961-0201-004-2)   
and   the   present   studies    were   supported   by   Basic  
Science 
Research Program, Ministry of Education, 1995 (BSRI-95-2413).

\end{document}